\documentclass[11pt,twoside,reqno]{amsart}

\usepackage{amsfonts,amsmath,amsthm}
\usepackage{latexsym}   

\theoremstyle{plain}
\newtheorem{thm}{Theorem}[section]

\newtheorem{prop}[thm]{Proposition}
\newtheorem{cor}[thm]{Corollary}
\newtheorem{lemma}[thm]{Lemma}

\theoremstyle{definition}
\newtheorem{definition}[thm]{Definition}

\theoremstyle{remark}
\newtheorem{example}[thm]{Example}
\newtheorem{rem}[thm]{Remark}

\numberwithin{equation}{section}

\def \N {{\mathbb N}}
\def \R {{\mathbb R}}

\def \eps {\varepsilon}
\def \d {\delta}

\def \t {\tau}

\newcommand{\dopu}{{:}\allowbreak\ }

\def \< {\langle}
\def \> {\rangle}

\newcommand{\loglike}[1]{\mathop{\rm #1}\nolimits}

\def \innt {\loglike {int}}

\newcommand{\DP}{Daugavet property}

\newcommand{\Lip}{\mathrm{Lip}}

\newcommand{\nnnorm}{|\mkern-2mu|\mkern-2mu|}

\newcommand{\ab}{\allowdisplaybreaks}
\newcommand{\begsta}{\begin{statements}}
\newcommand{\begaeq}{\begin{aequivalenz}}
\def\endsta{\end{statements}}
\def\endaeq{\end{aequivalenz}}
\newcommand{\bea}{\begin{eqnarray*}}
\newcommand{\eea}{\end{eqnarray*}}
\newcommand{\kref}[1]{(\ref{#1})}

\newcounter{abc}   
\newcounter{iiiii} 

\newenvironment{aequivalenz}
{\setcounter{iiiii}{0}
\begin{list}%
{{\rm (\roman{iiiii})}}
{\usecounter{iiiii}
\parsep=0pt plus 1pt
\topsep=1pt plus 2pt minus 1pt
\itemsep=1pt plus 2pt minus 1pt
\leftmargin=3\baselineskip
\labelsep=.6\baselineskip
\labelwidth=2.4\baselineskip
\rightmargin 0pt}%
}%
{\end{list}}

\newenvironment{statements}%
{\setcounter{abc}{0}
\begin{list}%
{{\rm (\alph{abc})}}
{\usecounter{abc}
\parsep=0pt plus 1pt
\topsep=1pt plus 2pt minus 1pt
\itemsep=1pt plus 2pt minus 1pt
\leftmargin=3\baselineskip
\labelsep=.6\baselineskip
\labelwidth=2.4\baselineskip
\rightmargin 0pt}%
}%
{\end{list}}

\begin{document}

\title{The Daugavet property for spaces of Lipschitz functions}

\author{Yevgen Ivakhno, Vladimir Kadets and Dirk Werner}

\address{Faculty of Mechanics and Mathematics, Kharkov National
University,\linebreak
 pl.~Svobody~4,  61077~Kharkov, Ukraine}
\email{ivakhnoj@yandex.ru}
\address{Faculty of Mechanics and Mathematics, Kharkov National
University,\linebreak
 pl.~Svobody~4,  61077~Kharkov, Ukraine}
\email{vova1kadets@yahoo.com}
\address{Department of Mathematics, Freie Universit\"at Berlin,
Arnimallee~2--6, \qquad {}\linebreak D-14\,195~Berlin, Germany}
\email{werner@math.fu-berlin.de}

\thanks{The work of the second-named author
was supported by a fellowship from the \textit{Alexander-von-Humboldt
Stiftung}.
}



\maketitle

\begin{abstract}
For a compact metric space $K$ the space $\Lip(K)$ has the Daugavet
property if and only if
the norm of every $f \in \Lip(K)$ is attained locally. If
$K$ is  a subset of an $L_p$-space, $1<p<\infty$, 
this is equivalent to the convexity of~$K$.
\end{abstract}

\thispagestyle{empty}

\section{Introduction}

A Banach space $X$ is said to have the \textit{Daugavet property}
if
\begin{equation}\label{DE}
\|\mathrm{Id} + T\| = 1 + \|T\|
\end{equation}
for every rank-1 operator $T\dopu X\to X$; then \kref{DE} also holds
for all weakly compact operators on $X$ and even all operators that do
not fix copies of $\ell_1$. The \DP\ was introduced in \cite{KSSW}
and further  studied in \cite{Shv1} and  \cite{KSW}, but examples of
spaces having the \DP\ have long been known; e.g., $C[0,1]$,
$L_1[0,1]$, $L_\infty[0,1]$, the disk algebra, $H^\infty$, etc.

In this paper we shall investigate the \DP\ for spaces of Lipschitz
functions. Throughout, $(K,\rho)$ stands for a complete metric space
that is not reduced to a singleton. The space of all Lipschitz
functions on $K$ will be equipped with the seminorm
$$
\|f\| = \sup\biggl\{ \frac{ |f(t_1)-f(t_2)| }{ \rho(t_1,t_2) }\dopu
t_1\neq t_2 \in K \biggr\}.
$$
If one quotients out the kernel of this seminorm, i.e., the constant
functions, one obtains the Banach space $\Lip(K)$, whose norm will
also be denoted by $\|\,\,.\,\,\|$. Equivalently, one can fix a point
$t_0\in K$ and consider the Banach space $\Lip_0(K)$ consisting of all
Lipschitz functions on $K$ that vanish at $t_0$, with the Lipschitz
constant as an actual norm. It is easily seen that $\Lip(K)$ and
$\Lip_0(K)$ are isometrically isomorphic.
In this paper we prefer the first point of
view, but will refer to the elements of $\Lip(K)$ as functions rather
than equivalence classes, as is familiar with $L_p$-spaces.

Since $\Lip[0,1]$ is isometric to $L_\infty[0,1]$ via differentiation
almost everywhere, it is clear that  $\Lip[0,1]$ has the \DP.
On the other hand the H\"older space $H^\alpha[0,1]$, being the dual
of a space with the RNP \cite[p.~83]{Weaver},
fails the \DP\ by the results of
\cite{Woj92}; $H^\alpha[0,1]$ is just the Lipschitz space for $K=[0,1]$
with the metric $\rho_\alpha(s,t)= |s-t|^\alpha$. But for the unit
square $Q=[0,1]  \times [0,1]$ with the Euclidean metric it is far from
obvious whether the \DP\ holds for $\Lip(Q)$; in fact, this will turn
out to be true as a special case of Theorem~\ref{thm1} below. The
validity of the \DP\ of $\Lip(Q)$ was asked in \cite{Dirk-IMB}.

Whereas for the ``classical'' function spaces the validity of the \DP\
is equivalent to a nonatomicity condition (\cite{FoiSin} for $C(K) $
and $L_1(\mu)$, \cite{Woj92} for function algebras, \cite{Dirk10} for
$L_1$-preduals and \cite{Oik} for the noncommutative case), in the
setting of Lipschitz spaces it is a locality condition that plays a
similar role, for in Theorem~\ref{thm2} we will show for a compact
metric space $K$ that the \DP\ of $\Lip(K)$ is equivalent to the fact
that every Lipschitz function on $K$ almost attains its norm at
close-by points; see Definition~\ref{local}(a) for precision.
We also characterise compact ``local'' metric spaces by a condition
that is reminiscent of metric convexity (Proposition~\ref{propZ2})
and is sometimes even equivalent to it, e.g., for compact subsets of
$L_p$, $1<p<\infty$  (Proposition~\ref{propZ3}). As a result, 
for a compact subset of $L_p$, $1<p<\infty$, the 
\DP\ of $\Lip(K)$ is equivalent to the convexity of $K$.

An important tool to construct Lipschitz functions is McShane's
extension theorem saying that if $M\subset K$ and $f\dopu M\to \R$ is
a Lipschitz function, then there is an extension to a Lipschitz function
$F\dopu K\to \R$ with the same Lipschitz constant; see
\cite[p.~12/13]{BL1}. This will be used several times.

We will also make use of the following geometric characterisations of
the \DP\ from \cite{KSSW} and \cite{BKSW}. Part~(iii) is particularly
useful when one doesn't have full access to the dual space.
As for notation, we denote the closed unit ball (resp.\ sphere)
of a Banach space $X$ by $B_X $ (resp.\ $S_X$) and the closed ball
with centre $t$ and radius $r$ in a metric space $K$ by $B_K(t,r)$.

\begin{lemma} \label{l:MAIN}
The following assertions are equivalent:
\begin{aequivalenz}
\item
     $X$ has the Daugavet property.
\item
     For every $y\in S_X$, 
$x^*\in S_{X^*}$  and $\eps >0$ there
         exists some $x\in S_X$  such that $x^*(x) \ge 1-\eps$  and
$\|x+y \| \ge 2-\eps$.
\item
For every $\eps >0$ and for every $y \in S_X$ the
closed convex hull of the set
$\{u \in (1+\eps)B_X\dopu  \|y+u\| \ge 2-\eps\}$
contains $S_X$.
\end{aequivalenz}
\end{lemma}

\section{Local metric spaces}
\label{sec2}

Let us recall that a metric space $K$ is called \textit{metrically
convex} if for any two points $t_1,t_2\in K$ two closed balls
$B_K(t_1,r_1)$ and $B_K(t_2,r_2)$ intersect if and only if
$\rho(t_1,t_2)\le r_1+r_2$.

Clearly, convex subsets of normed spaces are metrically convex, and
$S^1=\{(x,y)\in \R^2\dopu x^2+y^2=1\}$ is metrically convex for the
geodesic metric, but not for the Euclidean metric.

We shall need the following lemma.

\begin{lemma} \label{M-conv}
A complete metric space $K$ is
metrically convex if and only if  for every two distinct points $t, \t \in K$
there is an isometric embedding $\phi\dopu  [0, a] \to K$ $($where
$a=\rho(t, \t))$ such that $\phi(0) = t$, $\phi(a) = \t$. In other
words, $K$ is metrically convex if and only if
every two points of $K$ can be connected by an isometric
copy of a linear segment.
\end{lemma}

\begin{proof}
The property displayed in the lemma clearly implies the metric
convexity of~$K$. To prove the converse, let $K$ be metrically convex
and let $t$ and $\tau$ be two points at a distance~$a$; we shall label
them $t_0$ and $t_a$. Then there is a point
$t_{a/2} \in B_K(t_0, a/2) \cap B_K(t_a, a/2)$. It follows that
$\rho(t_0,t_{a/2})= \rho(t_{a/2}, t_a) = a/2$. Likewise, pick points
$t_{a/4} \in B_K(t_0, a/4) \cap B_K(t_{a/2}, a/4)$ and
$t_{3/4\cdot a} \in B_K(t_{a/2}, a/4) \cap B_K(t_{a}, a/4)$.
Continuing in this manner, one obtains for each dyadic rational $d\in
[0,1]$ a point $t_{da}\in K$ such that $\rho(t_ {da}, t_{d'a})=
|d-d'|a$. The mapping $da\mapsto t_{da}$ can now be extended to an
isometric mapping $\phi\dopu [0,a]\to K$, as requested.
\end{proof}

The following definition is crucial for this paper.

\begin{definition} \label{local}
Let $K$ be a metric space.
\begsta
\item
The space $K$ is called \textit{local} if for every $\eps > 0$ and for
every function $f \in \Lip(K)$ there are two distinct points
 $\t_1,\t_2\in K$ such that $\rho(\t_1,\t_2) < \eps$ and
\begin{equation} \label{eq6}
\frac{f(\t_2) - f(\t_1)}{\rho(\t_1,\t_2)} > \|f\| - \eps.
\end{equation}
\item
Let $f \in \Lip(K)$ and $\eps > 0$. A point $t \in K$ is said to be an
\textit{$\eps$-point} of $f$ if in every neighbourhood $U \subset K$ of $t$
there are two points $\t_1,\t_2 \in U$ for which (\ref{eq6}) holds
true.
\item
The space $K$ is called \textit{spreadingly local} if for every
$\eps > 0$ and for every function $f \in \Lip(K)$ there are
infinitely many $\eps$-points of~$f$.
\endsta
\end{definition}

The next proposition provides a large class of examples.

\begin{prop} \label{m-conv-loc}
A metrically convex complete metric space $K$ is spreadingly
local.
\end{prop}

\begin{proof}
Fix an $\eps > 0$ and a function $f \in \Lip(K)$ with $\|f\| = 1$.
Select $t, \t \in K$ with $\rho(t, \t) > 0$ such that
$$
f(\t) - f(t) > (1-\eps)\rho(t, \t).
$$
Denote $a=\rho(t, \t)$ and apply Lemma~\ref{M-conv} to
this pair of points. The function $F = f \circ \phi\dopu  [0, a] \to
\R$, where $\phi$ is from Lemma~\ref{M-conv},
is 1-Lipschitz. Hence $|F'| \le 1$ a.e.\  on $[0, a]$ and
$$
\int_0^a F'(r) \,dr = f(\t) - f(t) > (1 - \eps)a.
$$

Therefore there are infinitely many points $r_i \in [0, a]$ with $F'(r_i)
> 1-\eps$. Let us show that every point of the form
$t_i=\phi(r_i)$ is an $\eps$-point of $f$.  By the definition of the
derivative we have
$$
\frac{F(r_i + \delta_i) - F(r_i)}{\delta_i} > 1-\eps.
$$
for sufficiently small $\delta_i \in (0,\eps)$.
Denote $\t_i=\phi(r_i+ \delta_i)$. Then $\rho(t_i, \t_i)= \delta_i$ and
$f(\t_i) - f(t_i) > (1-\eps)\delta_i$.
\end{proof}

Actually this proposition applies to a slightly more general class of
spaces $K$, defined by the requirement that for each pair of points
$t,\tau\in K$ and each $\eta>0$ there exists a curve of length ${\le
  \rho(t,\tau)} +\eta =: a_\eta$ joining $t$ and $\tau$. In other
words, there exists a $1$-Lipschitz mapping (having arclength as
parameter) $\phi\dopu [0,a_\eta]\to
K$ with $\phi(0)=t$, $\phi(a_\eta)=\tau$. Such spaces could be termed
\textit{almost metrically convex}. A variant of the above proof then
shows that almost metrically convex spaces are spreadingly local.

\begin{example} \label{exZ}
There is a (noncompact) almost metrically convex space that is not
metrically convex. Indeed, let
$$
M=\{f\in L_1[0,1]\dopu |f|=1 \mbox{ a.e.}\};
$$
this is a closed subset of $L_1$. Instead of the $L_1$-norm we shall
use the following equivalent norm on $L_1$. Pick a total sequence of
functionals $x_n^*\in S_{L_\infty}$ and put, for $f\in L_1$,
$$
\nnnorm f \nnnorm = \|f\|_{L_1} + \biggl( \sum_{n=1}^\infty 2^{-n}
|x_n^*(f)|^2  \biggr)^{1/2}.
$$
This norm is strictly convex. It follows that $M$, equipped with the
metric $\rho(f,g )= \nnnorm f-g \nnnorm$, is not metrically convex
since it is not convex; indeed, if $f,g\in M$, then \textit{no}
nontrivial convex combination belongs to $M$ (unless $f=g$).

On the other hand, $(M,\rho)$ is almost metrically convex. To 
see this let $f\neq g$ be two functions in~$M$. 
For a Borel set $A\subset [0,1]$ define
$h_A\in M$ by 
$$
h_A = f\chi_A + g \chi_{[0,1]\setminus A}.
$$
Given $\eps>0$, pick $\eps' \le \eps \nnnorm f-g \nnnorm $ and
$N\in\N$ such that $2 ( \sum_{n>N} 2^{-n} )^{1/2} \le \eps'$.
Define a nonatomic vector measure taking values in $\R^{N+1}$ by
$$ \textstyle
\mu(A) = \bigl( \int_A |f-g|, x_1^*((f-g)\chi_A), \ldots,
x_N^*((f-g)\chi_A) \bigr).
$$
By the Lyapunov convexity theorem \cite[Th.~5.5]{RudFA} there exists a
Borel set $\Delta$ such that $\mu(\Delta)= \frac12 \mu([0,1])$. We then
have, since $g-h_\Delta = (g-f)\chi_\Delta$
\bea
\nnnorm g-h_\Delta \nnnorm &=&
\| g-h_\Delta \|_{L_1} + \biggl( \sum_{n=1}^\infty 2^{-n}
|x_n^*( g-h_\Delta )|^2  \biggr)^{1/2} \ab  \\
&\le&
\int_\Delta |f-g| + \biggl( \sum_{n=1}^N 2^{-n}
|x_n^*( (f-g)\chi_\Delta )|^2  \biggr)^{1/2} + \eps' \ab \\
&=&
\frac12 \int_0^1 |f-g| + \frac12 \biggl( \sum_{n=1}^N 2^{-n}
|x_n^*( f-g )|^2  \biggr)^{1/2} + \eps' \\
&\le&
\frac12 \nnnorm f-g \nnnorm + \eps' \le
\Bigl( \frac12 + \eps \Bigr) \nnnorm f-g \nnnorm 
\eea
and likewise
$$
\nnnorm f-h_\Delta \nnnorm \le
\Bigl( \frac12 + \eps \Bigr) \nnnorm f-g \nnnorm .
$$
Let $F_0=f$, $F_1=g$, $F_{1/2}=h_\Delta$. Now we reiterate the above
construction, first applying it to $F_0$, $F_{1/2}$ and $\eps/2$ and
then to $F_{1/2}$, $F_1$ and $\eps/2$ to obtain functions $F_{1/4},
F_{3/4}\in M$ such that
\bea
\max\{ \nnnorm F_0 - F_{1/4} \nnnorm, \nnnorm F_{1/2} - F_{1/4}
\nnnorm \} &\le &
\Bigl( \frac12 + \frac \eps2 \Bigr) \nnnorm F_0 - F_{1/2} \nnnorm, \\
\max\{ \nnnorm F_{1/2} - F_{3/4} \nnnorm, \nnnorm F_{1} - F_{1/4}
\nnnorm \} &\le &
\Bigl( \frac12 + \frac \eps2 \Bigr) \nnnorm F_1 - F_{1/2} \nnnorm.
\eea
Continuing in this manner, we can assign to each dyadic rational
$d\in[0,1]$ a function $F_d\in M$ such that the curve $[0,1]\to M$,
$t\mapsto F_t$, obtained from this by continuous extension, has a
length that can be estimated from above by 
$$
\sup_n \Bigl( \frac12 + \frac\eps{2^{n-1}} \Bigr)
\Bigl( \frac12 + \frac\eps{2^{n-2}} \Bigr)
\cdots
\Bigl( \frac12 + \eps \Bigr) 2^n
\le \exp (2^{2-n}\eps + \cdots + 2\eps)
\le e^{4\eps}.
$$
Therefore $M$ is almost metrically convex.
\end{example}

We will need a lemma in order to control the Lipschitz constant of a
function by the Lipschitz constant of some restriction under highly
technical assumptions that we shall meet later.
In the following, $\sqcup$ is used to indicate a disjoint union.

\begin{lemma}\label{flat}
Let $A=B \sqcup C$ be a metric space, $r \in (0,1/4]$, $\d <
r^2/16$, $\rho(B,C) > r$. Suppose $\tilde{C} \subset C$ is a $\d$-net of
$C$ such that every two points of $\tilde{C}$ are at least
$r$-distant, and let $f\dopu  A \to \R$ be  a function that  is
$1$-Lipschitz on $B \sqcup \tilde{C}$ and  also $1$-Lipschitz on
every ball $B_A(t, \d)$ for $t\in \tilde{C}$. Then $f$ is
$(1 + r/2)$-Lipschitz on the whole space~$A$.
\end{lemma}

\begin{proof}
Consider arbitrary points $s_1 \neq s_2 \in A$. We have to prove that
\begin{equation} \label{eq5}
\left|\frac{f(s_2) - f(s_1)}{\rho(s_1,s_2)}\right| \le 1 + r/2.
\end{equation}
We have to distinguish three cases:  firstly, when $s_1, s_2 \in B$;
secondly, when $s_1, s_2 \in C$; and  thirdly, when one of the
points (say, $s_1$) belongs to $B$ and the other belongs to $C$.

In the first case (\ref{eq5}) holds true even with 1 on the right
hand side by assumption on $f$.
Consider the second case. If $s_1, s_2$ belong to
the same ball of the form $B_A(t, \d)$ for $t\in \tilde{C}$, then
the job is likewise done. If not, let $t_1 \neq t_2 \in \tilde{C}$ be
points such that $\rho(t_1, s_1) \le \d$ and $\rho(t_2, s_2) \le  \d$.
Then
\bea
\left|\frac{f(s_2) - f(s_1)}{\rho(s_1,s_2)}\right|
&\le&
\left|\frac{f(s_2) - f(t_2)}{\rho(s_1,s_2)}\right| +
\left|\frac{f(t_2) - f(t_1)}{\rho(s_1,s_2)}\right|
+\left|\frac{f(t_1) - f(s_1)}{\rho(s_1,s_2)}\right| 
\ab  \\
&\le&
\frac{\d}{\rho(s_1,s_2)}
+\frac{\rho(t_2,t_1)}{\rho(s_1,s_2)}+ \frac{\d}{\rho(s_1,s_2)} 
\ab \\
&\le&
\frac{2\d}{r - 2\d} +\frac{\rho(t_2,t_1)}{\rho(t_2,t_1) - 2\d} \\
&\le&
\frac{2\d}{r - 2\d} +1 + \frac{2\delta}{\rho(t_2,t_1) - 2\d} \\
&\le&
1 + \frac{4\d}{r - 2\d} \le  1 + r/2.
\eea
In the last case find $t_2 \in \tilde{C}$ such that  $\rho(t_2,
s_2) \le \d$. Then
\bea
\left|\frac{f(s_2) - f(s_1)}{\rho(s_1,s_2)}\right|
&\le&
\left|\frac{f(s_2) - f(t_2)}{\rho(s_1,s_2)}\right| +
\left|\frac{f(t_2) - f(s_1)}{\rho(s_1,s_2)}\right| \\
&\le&
\frac{\d}{\rho(s_1,s_2)} +\frac{\rho(t_2,s_1)}{\rho(s_1,s_2)} \\
&\le&
\frac{\d}{r} +\frac{\rho(t_2,s_1)}{\rho(t_2,s_1) - \d } \\
&\le&
\frac{\d}{r} + \frac{r}{r - \d } =
1 + \frac{\d}{r}  + \frac{\d}{r - \d } \le 1 + r/2.
\eea
This completes the proof of the lemma.
\end{proof}

Obviously, a spreadingly local space is local. In the compact case the
converse is valid, too, as will be pointed out now.

\begin{lemma}\label{many e-points}
If $K$ is compact and  local, then it is spreadingly
local.
\end{lemma}

\begin{proof}
We will prove by induction on $n$ that for every $f \in \Lip(K)$ and
for every $\eps > 0$ there are $n$ $\eps$-points of $f$.

Thanks to the compactness of $K$ every function $f\in \Lip(K)$ has a
``0-point'', i.e., a point that is an $\eps$-point for every $\eps>0$.
Indeed, take a sequence of pairs $t_n,\t_n \in K$ satisfying
Definition~\ref{local} with $\eps = 1/n$, $n=1,2, \ldots\,$, and take
an arbitrary limit point of $(t_n)$. So the start of the induction
holds true. Now assume the statement for a fixed $n$ and let us
prove it for $n+1$.

Take an $f \in \Lip(K)$ with $\|f\| = 1$ and $\eps \in (0,1/4]$.  Due
to our hypothesis there are $\eps$-points $t_1,  \ldots , t_n
$ of $f$. Also, select two points $\t_1, \t_2 \in K$ distinct from
all the $t_i$ and such that
$$
\frac{f(\t_2) - f(\t_1)}{\rho(\t_1,\t_2)} > 1 - \eps/4.
$$
Let $r\in (0,\eps/4]$ be a number so small that the balls
$U_i=B_K(t_i,r)$, $i=1, \ldots, n$, are disjoint and
contain neither $\t_1$ nor $\t_2$. Fix a $\d < r^2/16$, denote
 the interior of $B_K(t_i, \d)$ by $V_i$ and consider $\tilde{K} = (K
\setminus \bigcup_{i=1}^n U_i) \sqcup \bigcup_{i=1}^n V_i$ as a
subspace of the metric space $K$. Define $\tilde{f}\dopu  \tilde{K} \to
\R$ as follows: $\tilde{f}(t) = f(t)$ for $t \in K \setminus
\bigcup_{i=1}^n U_i$ and $\tilde{f}(t) = f(t_i)$ on the
corresponding $V_i$.  Lemma~\ref{flat} implies that $\tilde{f}$
satisfies a Lipschitz condition on $\tilde{K}$ with the constant
$ {1 + \eps/2}$. Extend $\tilde{f}$ to a function on $K$
preserving the Lipschitz  constant, still denoted by $\tilde{f}$.

Take as $t_{n+1}$ an arbitrary 0-point of the function $g =
f+\tilde{f}$. Since
$$
\|g\| \ge \frac{g(\t_2) - g(\t_1)}{\rho(\t_1,\t_2)} =
2\frac{f(\t_2) - f(\t_1)}{\rho(\t_1,\t_2)} > 2 - \eps/2,
$$
in every neighbourhood of $t_{n+1}$ there are points $s_1, s_2$
with
\begin{equation} \label{eq4}
\frac{f(s_2) - f(s_1)}{\rho(s_1,s_2)} + \frac{\tilde{f}(s_2) -
\tilde{f}(s_1)}{\rho(s_1,s_2)} > 2 - \eps/2.
\end{equation}
This implies that $t_{n+1}$ cannot belong to any  $V_i$ since
in $V_i$ the second fraction of (\ref{eq4}) is zero, but the first
one is not greater than 1; hence $t_{n+1}$ differs from all the
other $t_i$. On the other hand, by our construction $\|\tilde{f}\|
\le 1 + \eps/2$, so the second fraction of (\ref{eq4}) is
$\le1 + \eps/2$. Hence there is an estimate for the first fraction,
namely
$$
\frac{f(s_2) - f(s_1)}{\rho(s_1,s_2)} > 1 - \eps,
$$
which means that $t_{n+1}$ is an $\eps$-point of $f$.
\end{proof}

Next we are going to  characterise local metric spaces
intrinsically, at least in the compact case, using the following
geometric property that we have chosen to give an ad-hoc name.

\begin{definition}\label{defiZ1}
A metric space $K$ has \textit{property $(Z)$} if the following
condition is met: Given $t,\tau\in K$ and $\eps>0$, there is
some $z\in K\setminus \{t,\tau\}$ satisfying 
\begin{equation} \label{eqZ1}
\rho(t,z) + \rho(z,\tau) \le \rho(t,\tau) + \eps
\min\{\rho(z,t),\rho(z,\tau)\}. 
\end{equation}
\end{definition}

A compact space satisfying \kref{eqZ1} with $\eps=0$ is easily seen to
be metrically convex. Thus, property~$(Z)$ is ``$\eps$-close'' to
metric convexity, and there are instances when $(Z)$ actually implies
metric convexity; see Corollary~\ref{corZ4} and Remark~\ref{remZ5}
below.

Here is the connection between locality and property~$(Z)$.

\begin{prop}\label{propZ2}
Let $K$ be a metric space.
\begsta
\item
If $K$ is local, then $K$ has property~$(Z)$.
\item
If $K$ is compact and has property~$(Z)$, then $K$ is local.
\endsta
\end{prop}

\begin{proof}
(a) Assume that $K$ fails property~$(Z)$, i.e., for some
$t_0,\tau_0\in K$ and $\eps_0>0$ there are no points $z\in K\setminus
\{t_0,\tau_0\}$ as in \kref{eqZ1}. For a point $z\in K$ let $r(z)=
\rho(z,t_0)$, $s(z)= \rho(z,\tau_0)$ and $d= \rho(t_0,\tau_0)$. Pick
$\eps>0$ with
$$
\frac \eps {1-\eps} < \frac{\eps_0}4 .
$$
Now define $f\dopu K\to \R$ by 
$$
f(z)= 
\left\{ 
\begin{array}{l}
\max\{d/2-(1-\eps)s(z),0\} \\
\mbox{\qquad if }r(z)\ge s(z), \ r(z) + (1-2\eps)s(z) \ge d, \\
-\max\{d/2-(1-\eps)r(z),0\}\\
\mbox{\qquad if }r(z)\le s(z), \ (1-2\eps)r(z) + s(z) \ge d.
\end{array}
\right.
$$ 
This function is well defined, since for $r(z)=s(z)$ both parts of the
definition yield~$0$, and all points of $K$ are covered in the two
``if'' cases by our assumption on $K$; note that $2\eps<\eps_0$.

Let us show that $f$ is a Lipschitz function with $\|f\|=1$. Indeed,
the only critical case is to estimate $f(z_2)-f(z_1)$ when $f(z_2)>0$
and $f(z_1)<0$; in this case
\bea
f(z_2)-f(z_1) 
&=&
\Bigl( \frac d2 - (1-\eps) s(z_2) \Bigr) +
\Bigl( \frac d2 - (1-\eps) r(z_1) \Bigr) \\
&\le&
\Bigl( \frac{ r(z_2)+(1-2\eps) s(z_2)  }2 - (1-\eps) s(z_2) \Bigr) \\
&&
\mbox{}\qquad + 
\Bigl( \frac{ (1-2\eps)r(z_1)+ s(z_1)  }2 - (1-\eps) r(z_1) \Bigr) \\
&=&
\frac 12 (r(z_2)-s(z_2)) + \frac12 (s(z_1)-r(z_1)) \\
&\le&
\rho(z_1,z_2);
\eea
also, the norm is attained at $\tau_0,t_0$, i.e., 
$f(\tau_0)-f(t_0) = \rho(\tau_0,t_0)$.

Consider now points $z_1,z_2\in K$ where
\begin{equation} \label{eqZ5}
\frac { f(z_2)-f(z_1) }{ \rho(z_2,z_1) } > 1-\eps;
\end{equation}
we shall show that then $z_1$ is close to $t_0$ and $z_2$ is close to
$\tau_0$ so that their distance is necessarily big. Obviously, we must
have $f(z_2)>0$ and $f(z_1)<0$ for \kref{eqZ5} to subsist. In
particular, we have
\begin{equation} \label{eqZ6}
\rho(z_1,t_0) < \rho(z_1,\tau_0); \quad \rho(z_2,\tau_0)<\rho(z_2,t_0) .
\end{equation}
Hence
\bea
(1-\eps)\rho(z_1,z_2) &<&
f(z_2)-f(z_1) \\
&=&
\Bigl( \frac d2 -(1-\eps) \rho(z_2,\tau_0) \Bigr) -
\Bigl( \frac d2 -(1-\eps) \rho(z_1,t_0) \Bigr)\\
&=&
d-(1-\eps)( \rho(z_2,\tau_0) + \rho(z_1,t_0));
\eea
in other words
$$
(1-\eps)( \rho(z_1,z_2) + \rho(z_2,\tau_0) + \rho(z_1,t_0) ) <d
$$
so that
\begin{equation} \label{eqZ7}
\rho(z_k,t_0) + \rho(z_k,\tau_0) < \frac d{1-\eps} , \quad k=1,2.
\end{equation}
By our choice of $\eps_0$, $t_0$, $\tau_0$ and \kref{eqZ6}
$$
\rho(z_1,t_0) + \rho(z_1,\tau_0) \ge d+\eps_0 \rho(z_1,t_0)
$$
so that by \kref{eqZ7}
$$
d+\eps_0 \rho(z_1,t_0) < \frac d{1-\eps} 
$$
and hence $\rho(z_1,t_0) < d/4$ by our choice of $\eps$. Likewise
$\rho(z_2,\tau_0) < d/4$ and consequently $\rho(z_1,z_2) > d/2$.
Therefore, $K$ cannot be local. 

(b) Assume that $K$ is not local. Then there is a Lipschitz function
$f$ with $\|f\|=1$ for which \kref{eq6} is impossible for
$\tau_1,\tau_2$ at small distance, viz.\ for
$\rho(\tau_1,\tau_2)<\eps$.
By a compactness argument one hence deduces the existence of points
$t,\tau\in K$ such that 
\begin{equation} \label{eqZ2}
\frac{ f(\tau)-f(t) }{ \rho(\tau,t) } =1 
\end{equation}
and $\rho(t,\tau)$ is minimal among all points as in \kref{eqZ2}. 
Now let $\eps_n \searrow 0$ and apply condition~$(Z)$ to $t,\tau$ and
$\eps_n$. This yields a sequence of points $z_n\in K\setminus
\{t,\tau\}$ such that 
\begin{equation} \label{eqZ2a}
\rho(t,z_n) + \rho(z_n,\tau) \le \rho(t,\tau) + \eps_n
\min\{\rho(z_n,t),\rho(z_n,\tau)\}. 
\end{equation}
Passing to a subsequence we may assume that $(z_n)$ converges, say
$z_n\to z_0$, and that without loss of generality 
\begin{equation} \label{eqZ3}
\rho(t,z_n) \le \rho(\tau,z_n) \qquad \forall n\ge1.
\end{equation}
Note that 
\begin{equation} \label{eqZ4}
\rho(t,z_0) + \rho(z_0,\tau) = \rho(t,\tau).
\end{equation}
If $z_0\neq t$, then
\bea
1\ge \frac { f(z_0)-f(t) }{ \rho(z_0,t) }
&=&
\frac{ f(\tau)-f(t) }{ \rho(\tau,t) } \, 
\frac{ \rho(\tau,t) }{ \rho(z_0,t) } -
\frac{ f(\tau)-f(z_0) }{ \rho(\tau,z_0) } \, 
\frac{ \rho(z_0,\tau) }{ \rho(z_0,t) } \\
&\ge&
\frac{ \rho(\tau,t) }{ \rho(z_0,t) } -
\frac{ \rho(z_0,\tau) }{ \rho(z_0,t) } 
=1
\eea
by \kref{eqZ4}, and thus $f$ attains its norm at the pair $z_0,t$. 
But by \kref{eqZ3}
$$
\rho(t,z_0) \le \frac12 ( \rho(t,z_0) + \rho(\tau,z_0) ) =
\frac12 \rho(t,\tau),
$$
which contradicts the minimality condition imposed on the pair
$t,\tau$.

Therefore, $z_n\to t$, and for sufficiently large $n$ we have
$\rho(t,z_n)<\eps $ along with \kref{eqZ2a}. But then
\bea
\frac { f(z_n)-f(t) }{ \rho(z_n,t) } 
&=&
\frac{ f(\tau)-f(t) }{ \rho(\tau,t) }\, 
\frac{ \rho(\tau,t) }{ \rho(t,z_n) } -
\frac{ f(\tau)-f(z_n) }{ \rho(\tau,z_n) } \, 
\frac{ \rho(\tau,z_n) }{ \rho(t,z_n) } \\
&\ge&
\frac{ \rho(\tau,t) - \rho(\tau,z_n) }{ \rho(t,z_n) } 
\ge 1-\eps
\eea
by \kref{eqZ2a}, which contradicts our choice of $f$, since
$\rho(t,z_n)<\eps $.
\end{proof}

The definition of locality immediately implies that a compact local
space is connected; one just has to apply the definition with the
indicator function of a set that is both open and closed.
We will now present a class of compact metric spaces for which
property~$(Z)$ and hence locality implies (metric) convexity. 
Recall that a Banach space $(E,\|\,\,.\,\,\|_E)$ is called
\textit{locally uniformly rotund} if for each $x\in S_E$ and $\eta>0$
there is some $\delta=\delta_x(\eta)>0$ such that $\|x-y\|_E\le\eta$
whenever $y\in B_E$ and $\|\frac12(x+y)\|_E \ge 1-\delta$.

\begin{prop}\label{propZ3}
Let $(E,\|\,\,.\,\,\|_E)$  be a smooth locally uniformly rotund Banach 
space and let $K\subset E$ be a compact subset with
property~$(Z)$. Then $K$ is convex.
\end{prop}

\begin{proof}
By a result of Vlasov (\cite{Vlasov}, \cite[Th.~2.2, p.~368]{Singer}) 
a compact Chebyshev subset of a
smooth Banach space is convex. If we assume that $K$ is not convex,
this means that there are two points $P,Q\in K$ and a ball $B$ whose
interior does not intersect $K$ with $P,Q\in \partial B$; we may
assume that $B$ is centred at the origin, $B=B_E(0,\alpha)$, and by
scaling that $\|P-Q\|_E=1$. 
Applying condition~$(Z)$ to $P,Q$ and an arbitrary $\eps>0$ yields
some $z=z(\eps)\in K\setminus \{P,Q\}$ as in \kref{eqZ1}. We may as
well assume that $z_0= \lim_{\eps\to0} z(\eps) $ exists; $z_0$ lies on
the line segment $[P,Q]$ by strict convexity of~$E$. Thus $z_0=P$ or
$z_0=Q$; without loss of generality let us assume the latter. Fix, for
the time being, $\eps$ and $z=z(\eps)$ and put $r=\|z-Q\|_E$
${(<1/2)}$.

Now consider $Q(\lambda)= \lambda P + (1-\lambda) Q$, $0\le \lambda
\le1$. Let us estimate $\|z-Q(\lambda)\|_E$ in order to derive a
contradiction. On the one hand we have, since $z\in K$ and thus
$\|z\|_E\ge\alpha$,
$$
\|z-Q(\lambda)\|_E \ge \|z\|_E - \|Q(\lambda)\|_E 
\ge \alpha - \|Q(\lambda)\|_E =: \varphi(\lambda).
$$
Now $\varphi $ is  a concave function of $\lambda$ with $\varphi(0)=0$
and 
$$
\varphi(1/2) = \alpha - \Bigl\| \frac12 (P+Q) \Bigr\| 
>0
$$
by strict convexity. Hence with $\sigma = 2 \varphi(1/2)$ 
\begin{equation} \label{eqZ8}
\|z-Q(r)\|_E \ge\varphi(r) \ge \sigma r.
\end{equation}
On the other hand, \kref{eqZ1} means that $z\in B_E(P, 1-r+\eps r)$;
therefore the point $w= \frac12 (z+Q(r))$ also belongs to this ball,
but $w\notin \innt B_E(Q, r-\eps r)$. In other words,
\begin{equation} \label{eqZ9}
\Bigl\| \frac{ (Q-z)+(Q-Q(r)) }2 \Bigr\|_E =
\Bigl\| Q - \frac{ z+Q(r) }2 \Bigr\|_E \ge r-\eps r.
\end{equation}

Specifically, let $\eta=\sigma/2$ and $0<\eps<\delta_{P-Q}(\eta)$.
Then \kref{eqZ9} and local uniform rotundity (note that 
$(Q-z)/r, (Q-Q(r))/r\in B_E$) imply that 
$$
\|z-Q(r)\|_E \le r\eta < r\sigma
$$
contradicting \kref{eqZ8}.
\end{proof}

Proposition~\ref{propZ3} applies in particular to $L_p$-spaces for
$1<p<\infty$ and most particularly to Hilbert spaces.

We can sum up the previous results as follows.

\begin{cor}\label{corZ4}
Let $K$ be a compact metric space. Then the following are equivalent:
\begsta
\item[\rm(1)]
$K$ is local;
\item[\rm(2)]
$K$ is spreadingly local;
\item[\rm(3)]
$K$ has property~$(Z)$.
\endsta
If $K$ is a subset of a smooth locally uniformly rotund Banach space,
then a further equivalent condition is:
\begsta
\item[\rm(4)]
$K$ is convex.
\endsta
\end{cor}

Another link between locality and metric convexity is provided by the
following technical remark. 

\begin{rem}\label{remZ5}
Let us say that $K$ satisfies $(Z')$ if in addition
to \kref{eqZ1} in Definition~\ref{defiZ1} we require that
$$
\rho(z,\tau) \le \rho(z,t).
$$
Since one can exchange the roles of $t$ and $\tau$ here, this means
that there is
one point as in \kref{eqZ1} that is closer to $\tau$ than to $t$ and
another one that is closer to $t$ than to $\tau$. It is then possible
to show that  $(Z')$ implies metric convexity
for compact spaces; see below. Hence
locality implies metric convexity for those compact metric spaces that
are symmetric in the sense that for any two points in $K$ there is an
isometry on $K$ swapping these two points.

To prove this remark, we rephrase property~$(Z')$ by saying that for
every $\eps>0$ and every $t,\tau \in K$ there exists some $z\in
K\setminus\{\tau\} $ such that
\begin{eqnarray}
\label{eqZ10}
(1-\eps)\rho(\tau,z) + \rho(t,z) &\le& \rho(t,\tau), \\
\label{eqZ11}
\rho(\tau,z) &\le& \rho(t,z).
\end{eqnarray}
The strategy of the proof will be to infer from this in the compact
case that for
every $\eps>0$ and every $t,\tau \in K$ there exists some $z\in K$
for which \kref{eqZ10} holds and
\begin{equation}\label{eqZ12}
\frac1{10} \rho(t,\tau) \le \rho(\tau,z) \le \frac9{10} \rho(t,\tau).
\end{equation}
If we let $\eps\to0$ and consider a limit point $z_0$ of the
$z=z(\eps)$ satisfying \kref{eqZ10} and \kref{eqZ12}, then we can be
certain that $z_0\neq t$ and $z_0\neq \tau$, but 
\begin{equation}\label{eqZ12a}
\rho(t,z_0) + \rho(z_0,\tau) = \rho(t,\tau).
\end{equation}
As remarked earlier this implies the metric convexity of the compact
space~$K$.

Let us now come to the details. Fix $t$, $\tau$ and $\eps$; we may
suppose that $\rho(t,\tau)=1$. Assume for a contradiction that we
cannot achieve \kref{eqZ10} and \kref{eqZ12} simultaneously.
Let 
$$
K_0= \{z\in K\dopu \mbox{\kref{eqZ10} and \kref{eqZ11} hold}\}.
$$
Since $K_0\neq\{\tau\}$ by property~$(Z')$, there is some $u\in
K_0$ such that $\rho(u,t)<1$, and therefore $\alpha :=
\min\{\rho(z,t)\dopu z\in K_0\}$ is attained at some $u_0\in
K_0\setminus \{\tau\}$. Then 
$(1-\eps)\rho(\tau,u_0) + \rho(u_0,t)\le1$ by \kref{eqZ10}.
Now define $0\le \tilde\eps \le\eps$ by 
\begin{equation}\label{eqZ13}
(1-\tilde\eps)\rho(\tau,u_0) + \rho(u_0,t) =1.
\end{equation}
If $\tilde\eps=0$, we have already found a point as in  \kref{eqZ12a},
and we are done. So we assume that $\tilde\eps>0$ in the sequel. Then
we can apply \kref{eqZ10} and \kref{eqZ11}, i.e., property~$(Z)$, with
$t$, $u_0$ and $\tilde\eps$ in place of $t$, $\tau$ and $\eps$. This
yields some $\tilde z\in K\setminus \{u_0\}$ with
\begin{eqnarray}
\label{eqZ14}
(1-\tilde\eps)\rho(u_0,\tilde z) + \rho(t,\tilde z) &\le& \rho(t,u_0), \\
\label{eqZ15}
\rho(u_0, \tilde z) &\le& \rho(t,\tilde z).
\end{eqnarray}
Next, add \kref{eqZ13} and \kref{eqZ14} to obtain
\begin{equation}\label{eqZ16}
(1-\tilde\eps) (\rho(\tau,u_0) + \rho(u_0,\tilde z) ) + 
\rho(t,\tilde z)  \le1.
\end{equation}
But $\rho(t,\tilde z) < \rho(t,u_0) =\alpha$, since $\tilde z \neq
u_0$ in \kref{eqZ14}; hence $\tilde z \notin K_0$. Now the previous
inequality, \kref{eqZ16} and $\tilde\eps\le \eps$ show that $\tilde z$
satisfies \kref{eqZ10}; therefore it must fail \kref{eqZ11}, i.e., 
\begin{equation}\label{eqZ166}
\rho(\tau, \tilde z) > \rho(t,\tilde z).
\end{equation}   
Also, recall that $u_0$ satisfies \kref{eqZ10} and that we have
assumed that \kref{eqZ10} and \kref{eqZ12} do not hold
simultaneously. This implies that 
$$
\rho(\tau,u_0) < 1/10 \mbox{ \ or \ } \rho(\tau,u_0) > 9/10
$$
and
$$
\rho(\tau,\tilde z) < 1/10 \mbox{ \ or \ } \rho(\tau,\tilde z) > 9/10.
$$
If $\rho(\tau,u_0) > 9/10$, then $\rho(t,u_0) > 9/10$ by \kref{eqZ11};
recall that $u_0\in K_0$. Then \kref{eqZ13} furnishes the
contradiction 
$$
1= (1-\tilde\eps)\rho(\tau,u_0) + \rho(u_0,t) >
(2-\tilde\eps)\frac9{10} > 1
$$
if, say, $\eps\le1/4$. The conclusion at this point is 
\begin{equation}\label{eqZ17}
\rho(\tau,u_0) < 1/10 .
\end{equation}   
On the other hand, if $\rho(\tau,\tilde z) < 1/10$, then 
$\rho(t,\tilde z) > 9/10$ by the triangle inequality, which
contradicts \kref{eqZ166}. Consequently
\begin{equation}\label{eqZ18}
\rho(\tau,\tilde z) > 9/10.
\end{equation}  
If we now use that $\tilde z$ satisfies \kref{eqZ14} and \kref{eqZ15},
we derive, for $\eps\le1/4$, that
$$
\rho(u_0,\tilde z) \le \rho(t,\tilde z) 
\le 1 - (1-\eps) \rho(\tau,\tilde z) \le\frac{13}{40}
$$
and hence the contradiction 
$$
\rho(\tau,t) \le \rho(\tau,u_0) + \rho(u_0,\tilde z) + \rho(\tilde z,t) 
<1.
$$
This completes the proof of  the remark.
\end{rem}

We do not know any example of a compact space with $(Z)$ that is not
metrically convex.

\section{Locality and the \DP}
\label{sec3}

We can now prove a sufficient criterion for $\Lip(K)$ to have the
\DP. In particular it turns out that for closed convex subsets of Banach
spaces $\Lip(K)$ has the \DP.

\begin{thm} \label{thm1}
If $K$ is  a spreadingly local metric space (in particular if $K$ is
a metrically convex  or a compact local metric space),
then $\Lip(K)$ has the \DP.
\end{thm}

\begin{proof}
For short write $X=\Lip(K)$.
Due to Lemma~\ref{l:MAIN}
it is sufficient to prove that
for every $\eps \in (0,1/4]$, and for every $f,g \in S_X$ the
closed convex hull of the set
$W = \{u \in (1+\eps)B_X\dopu  \|f+u\| \ge 2-\eps\}$
contains $g$.

In order to do this fix an $n \in \N$ and select $\eps/2$-points
$s_1,  \ldots , s_n $ of~$f$. Let $r\in (0,\eps/4]$ be a
 number so small that the balls $U_i=B_K(s_i,r)$, $i=1, \ldots, n$, are
disjoint. Fix a $\d < r^2/16$, and select  $t_i, \t_i \in B_K(s_i,
\d)$ such that
\begin{equation} \label{eq7}
 f(\t_i) - f(t_i) > (1-\eps/2)\rho(t_i, \t_i).
\end{equation}
Consider $K_i = (K \setminus  U_i) \sqcup  \{t_i,\tau_i\}$ as a subspace
of the metric space~$K$.  Define $u_i\dopu  K_i \to \R$ as follows:
$u_i(t_i) = g(t_i)$, $u_i(\t_i) = g(t_i)+ f(\t_i) - f(t_i)$ and
$u_i(s) = g(s)$ on the rest of $K_i$. It follows from  Lemma~\ref{flat}
that $u_i$ satisfies a Lipschitz condition on $K_i$ with
the constant $1 + r/2 < 1+\eps/2$. Extend $u_i$ to a function on
 $K$ preserving the Lipschitz  constant, still denoted by $u_i$.

Note that each  $u_i$ belongs to $W$. In fact $\|u_i\| \le
1+\eps$ by construction and
$$
\|f+u_i\| \ge \frac{(f+u_i)(\t_i) - (f+u_i)(t_i)}{\rho(\t_i,t_i)}
= 2\frac{f(\t_i) - f(t_i)}{\rho(\t_i,t_i)}
> 2-\eps.
$$

On the other hand the arithmetic mean of the $u_i$
(the simplest convex combination) approximates g, for
$$
\biggl\|g - \frac1n \sum_{i=1}^n u_i\biggr\| =
\frac1n \biggl\|\sum_{i=1}^n (u_i - g) \biggr\| \le
\frac{4 + 2\eps}{n}.
$$
The last inequality follows from the fact that each  $u_i - g$
has  norm  $\le \|u_i\|+\|g\| \le 2 + \eps$ and their supports $U_i$ are
disjoint.
\end{proof}

Finally we  address the question in how far our locality
conditions are necessary for the \DP; for compact spaces, this will
turn out to be the case (Theorem~\ref{thm2} below). The bulk of the
technical work will be done in the following lemma.

\begin{lemma} \label{basic1}
Suppose $\Lip(K) $ has the \DP. Then for every $t_1,t_2 \in K$ with
$\rho(t_1,t_2) = a > 0$, 
for every $f \in S_{\Lip(K)}$ with $f(t_2) -
f(t_1) = a$  $($i.e., $f$ attains its norm at the pair
$t_1,t_2)$ and for every $\eps > 0$ there are $\t_1 = \t_1(\eps)
, \t_2= \t_2(\eps) \in K$ with the following properties:
\begin{enumerate}
\item \label{pr1} $f(\t_2) - f(\t_1) \ge (1 - \eps)\rho(\t_1,\t_2)$;
\item \label{pr3}
$\rho(t_1,\t_2) - \rho(t_1,\t_1) \ge (1 - \eps)\rho(\t_1,\t_2)$,\\
$\rho(t_2,\t_1) - \rho(t_2,\t_2) \ge (1 - \eps)\rho(\t_1,\t_2)$;
\item \label{pr4} $\rho(\t_1,\t_2) \to 0$ as $\eps \to 0$.
\end{enumerate}
\end{lemma}

\begin{proof}
We shall abbreviate $\Lip(K)$ by $X$.
Consider the following functions $y_i \in X$:
$$
y_1 = f, \ y_2(t) = \rho(t_1,t), \ y_3(t) = -\rho(t_2,t).
$$
For all these functions we have
\begin{equation} \label{eq2.2}
y_i(t_2) - y_i(t_1) = a, \ \|y_i\| = 1.
\end{equation}
Then the arithmetic mean
$y = (y_1 + y_2 + y_3 )/3$ is of norm~1 as well. Consider
$x^* \in X^*$, with the action
\begin{equation} \label{eq2.2a}
x^*(x) = \frac{1}{a} (x(t_2) - x(t_1)).
\end{equation}
Clearly $\|x^*\|=1$.
Due to the \DP\  of $X$ there is, by Lemma~\ref{l:MAIN},
an $x \in S_X$ such that $x^*(x) >
1 - \eps$, i.e.,
\begin{equation} \label{eq2.5}
x(t_2) - x(t_1) > (1 - \eps)a,
\end{equation}
and at the same time $\|x - y\| > 2 - \eps/3$. The last condition
means that there are two distinct points $\t_1,\t_2 \in K$ for
which
 $$
(x - y)(\t_1) - (x - y)(\t_2) > (2 - \eps/3)\rho(\t_1,\t_2),
$$
i.e.,
$$
\frac13 \sum_{i=1}^3 \left((x - y_i)(\t_1) - (x -
y_i)(\t_2)\right)> (2 - \eps/3)\rho(\t_1,\t_2).
$$
Since neither of these three summands  exceeds 2$\rho(\t_1,\t_2)$,
we get the following three inequalities:
\begin{equation} \label{eq2.3}
(x - y_i)(\t_1) - (x - y_i)(\t_2)> (2 - \eps)\rho(\t_1,\t_2),
\quad i=1,2,3.
\end{equation}
Taking into account $x(\t_1) - x (\t_2) \le \rho(\t_1,\t_2)$ we
deduce that
\begin{equation} \label{eq2.4}
 y_i(\t_2) -  y_i(\t_1) > (1 - \eps)\rho(\t_1,\t_2),
\quad i=1,2,3.
\end{equation}

The case $i=1$ gives us the requested property~\kref{pr1}, and 
the cases $i=2, 3$ of (\ref{eq2.4}) immediately provide
property~\kref{pr3}.
Finally, substituting the Lipschitz conditions $x(\t_1) \le x(t_1) +
\rho(t_1,\t_1)$ and  $x(\t_2) \ge x(t_2) - \rho(t_2,\t_2)$ into
(\ref{eq2.3}) and applying (\ref{eq2.5}) we obtain
\bea
(2 - \eps)\rho(\t_1,\t_2)
&<&
x(t_1) - x(t_2) + \rho(t_1,\t_1) + \rho(t_2,\t_2) + y_i(\t_2) -
y_i(\t_1) \\
&\le&
-(1 - \eps)\rho(t_1,t_2) + \rho(t_1,\t_1) + \rho(t_2,\t_2) + \rho(\t_1,\t_2),
\eea
so
\bea
(1 - \eps)\rho(t_1,t_2)
&<&
\rho(t_1,\t_1) + \rho(t_2,\t_2) - (1 - \eps)\rho(\t_1,\t_2)\\
&\le&
 (2 - \eps)\left(\rho(t_1,\t_1) +
\rho(t_2,\t_2)\right) -(1 - \eps)\rho(t_1,t_2)
\eea
by the triangle inequality; hence
$$
2\rho(t_1,\t_1) + 2\rho(t_2,\t_2) >
4(1 - \eps)/(2 - \eps)\rho(t_1,t_2).
$$
 Adding to this inequality both  inequalities from
property~\kref{pr3} we obtain
\bea
\lefteqn{ \mbox{\hspace*{-1cm}}
\rho(t_1,\t_1) + \rho(t_2,\t_2) + \rho(t_1,\t_2) + \rho(t_2,\t_1)
} \\
&\ge &
4(1 - \eps)/(2 - \eps)\rho(t_1,t_2) + 2(1 - \eps)\rho(\t_1,\t_2).
\eea
Since the left hand side is not greater than $2\rho(t_1,t_2)$ we deduce
$$
2(1 - \eps)\rho(\t_1,\t_2) \le
\biggl(2 - 4\frac{1 - \eps}{2 - \eps} \biggr)\rho(t_1,t_2)
$$
which gives  property~\kref{pr4}.
\end{proof}

We can now deduce
the   main theorem of this paper.

\begin{thm} \label{thm2}
If $K$ is a compact metric space, then $\Lip(K)$ has the \DP\  if and
only if $K$ is local.
\end{thm}

\begin{proof}
The ``if'' part has already been proved in  Theorem~\ref{thm1}. Let
us prove the ``only if'' part. Assume $K$ is not local. Then
there is a function $f \in \Lip(K)$, $\|f\|=1$, and there is an $r >
0$ such that
\begin{equation} \label{eq.neu}
f(\t_2) - f(\t_1) < (1 - r)\rho(\t_1,\t_2)
\end{equation}
for
every $\t_1, \t_2 \in K$ with $\rho(\t_1, \t_2) < r$. Hence by a
compactness argument there is a pair of points $t_1,t_2 \in K$
with $\rho(t_1,t_2)  > 0$ on which $f$ attains its norm, i.e.,
with $f(t_2) - f(t_1) = \rho(t_1,t_2)$. If nevertheless $\Lip(K)$ has the
\DP, then applying Lemma~\ref{basic1} to $f$ and these $t_1,t_2$
with $\eps \to 0$ entails a contradiction between \kref{eq.neu}
and properties \kref{pr1} and \kref{pr4} from the
lemma.
\end{proof}

The space $\Lip(K)$ has a canonical predual, called the Arens-Eells
space in \cite{Weaver} and the Lipschitz free space in \cite{GK-Lip}
and \cite{Kal-Lip}. Since we have used in \kref{eq2.2a}, in the proof of
Lemma~\ref{basic1}, a functional
from that predual, i.e., a weak$^*$ open slice, the lemma works under the
assumption that the Lipschitz free space on $K$ has the
\DP. Consequently, for a compact metric space $\Lip(K)$ has the \DP\
if and only if its Lipschitz free space has.

In the setting of subsets of certain Banach spaces like $L_p$,
 $1<p<\infty$,
 we can rephrase Theorem~\ref{thm2}
as follows, using Corollary~\ref{corZ4}.

\begin{cor} \label{corthm3}
If $K$ is a compact subset of a smooth locally uniformly rotund Banach space, 
then $\Lip(K)$ has
the \DP\  if and only if $K$ is convex.
\end{cor}


\end{document}